\documentclass[12pt]{amsart}
\input xypic
\xyoption{all}
\usepackage{yfonts}
\usepackage{graphicx}

\usepackage{amssymb}
\newtheorem{thm}{Theorem}[section]
\newtheorem{lemma}[thm]{Lemma}
\newtheorem{cor}[thm]{Corollary}
\newtheorem{prop}[thm]{Proposition}

\theoremstyle{definition}
\newtheorem{rem}[thm]{Remark}

\def\can{\operatorname{can}\nolimits}

\def\Comp{\operatorname{Comp}\nolimits}
\def\End{\operatorname{End}\nolimits}

\def\mgr{{\operatorname{\!-gr}\nolimits}}
\def\grdim{\operatorname{grdim}\nolimits}
\def\Ho{\operatorname{Ho}\nolimits}
\def\Hom{\operatorname{Hom}\nolimits}
\def\id{\operatorname{id}\nolimits}
\def\Id{\operatorname{Id}\nolimits}

\def\opp{{\operatorname{opp}\nolimits}}
\def\rev{{\operatorname{rev}\nolimits}}
\def\mMod{\operatorname{\!-mod}\nolimits}
\def\dotimes{\otimes\hspace{-0.485cm}\bigcirc}

\def\Gsl{{\mathfrak{sl}}}

\def\CA{{\mathcal{A}}}
\def\CB{{\mathcal{B}}}
\def\CC{{\mathcal{C}}}

\def\CL{{\mathcal{L}}}
\def\CM{{\mathcal{M}}}

\def\CT{{\mathcal{T}}}
\def\CU{{\mathcal{U}}}
\def\CV{{\mathcal{V}}}
\def\CW{{\mathcal{W}}}

\def\BQ{{\mathbf{Q}}}

\def\BZ{{\mathbf{Z}}}
\def\GS{{\mathfrak{S}}}

\def\tR{{\tilde{R}}}

\def\iso{\buildrel \sim\over\to}

\begin{document}

\author{Mark Ebert}
\address{M.E.: UCLA Mathematics Department, Los Angeles, CA 90095-1555, USA.}
\email{markebert@math.ucla.edu}
\author{Rapha\"el Rouquier}
\address{R.R.: UCLA Mathematics Department, Los Angeles, CA 90095-1555, USA.}
\email{rouquier@math.ucla.edu}
\thanks{The authors gratefully acknowledge support from the Simons Foundation collaboration
grant on New Structures in Low-dimensional Topology (grant \#994340).
The second author also gratefully acknowledges support from the NSF
(grant DMS-2302147) and from the Simons Foundation (grant \#376202).}
%\dedicatory{Preliminary version}

\title{Higher tensor product for $\Gsl_2$ and Webster algebras}

\date{\today}

\begin{abstract}
	We construct a model for the tensor product of the regular $2$-representation of
	the enveloping algebra of $\Gsl_2^+$ with the vector $2$-representation, based on the
	$\infty$-categorical definition of \cite{Rou3}. Our model contains McMillan's minimal one
	\cite{Mc}. Our use of an infinite family of generators provides a simpler model that we
	prove is equivalent to Webster's tensor product category \cite{We}.
\end{abstract}

\maketitle
\setcounter{tocdepth}{3}
\tableofcontents

\section{Introduction}

Higher representation theory is a version of representation theory where vector spaces are
replaced by categories. In this article, we consider the case of $\Gsl_2$, where the original
theory was introduced in \cite{ChRou}, and the graded version we consider in \cite{Lau}.

In \cite{Rou3}, the second author defines a tensor product for $2$-representations and conjectures
that the tensor product of simple $2$-representations agrees with
Webster's quiver Hecke algebra $2$-representations.

\smallskip
The tensor product of $2$-representations involves only the positive part of the Lie algebra and
in this article we consider the case of the one-dimensional Lie algebra $\Gsl_2^+$. A 
$2$-representation of $\Gsl_2^+$ on a category over a field $k$
is the data of an endofunctor $E$ and of
endomorphisms $x$ of $E$ and $\tau$ of $E^2$ satisfying Hecke relations (cf (\ref{eq:xtau})
in \S\ref{se:monoidalcat}).
Equivalently, it
is the action of the monoidal category $\CU$ generated by $E$ and endomorphisms $x$ and $\tau$.
The endomorphism algebras in $\CU$ are nil affine Hecke algebras of type $A$. Note that
$\CU$ is a $2$-representation by left tensor product. The $2$-dimensional vector representation has a
categorical version, $\CL=\CL_0\oplus\CL_1$, with $\CL_i$ the category with one object whose
endomorphism ring is $k[y]$.

\smallskip
We consider the tensor product $\CU\dotimes\CL$ of $\CU$ and $\CL$.
The general theory defines this as an
$\infty$-category with an $\infty$-categorical $2$-representation. It also provides
a dg-model for the category together with an endofunctor $E$, and that is our starting point.
The endomorphisms $x$ and $\tau$ exist only on the derived category. We construct
a new $t$-structure on the derived category for which the action of $E$ is exact and we use
the compatibility $\CU\dotimes\CL\twoheadrightarrow \CU\dotimes\CL_0=\CU$
arising from the categorical version of the one-dimensional quotient of the vector
representation to determine $x$ and $\tau$, in the new $t$-structure.

This $t$-structure was considered by McMillan \cite{Mc} who provides an intricate
description of the 
endomorphism ring of a small progenerator in the heart and the bimodule corresponding to $E$, as well
as the endomorphisms $x$ and $\tau$.
What we do here instead is consider an infinite progenerating family for this $t$-structure.
We show that the full subcategory with those as objects
has a very simple description: it is one of Webster's tensor product categories.

\smallskip
Our work is a step in the description of the braided monoidal category of $2$-representations
of $\Gsl_2$, toward fulfilling Crane and Frenkel's program \cite{CrFr}.

\section{Notations}
Let $\phi:A\to B$ be a morphism of rings. Given $M$ a left (resp.
right) $B$-module, we denote by $\phi M$ (resp. $M\phi$) the left (resp. right) $A$-module
whose underlying set is $M$ and where the action of $a\in A$ is that of $\phi(a)$.

We denote by $A\mMod$ the category of finitely generated $A$-modules and by
$D^b(A)$ its derived category (when $A\mMod$ is abelian).

\smallskip
We fix a field $k$. By category, we mean a category enriched in $k$-vector spaces and
all functors between such categories are assumed to be $k$-linear.

Given a category $\CC$, we denote by $\CC^i$
the full subcategory of the category of functors from $\CC^\opp$ to the category of
sets with objects the direct summands of finite direct sums of objects of the form
$\Hom(-,c)$ for $c\in\CC$. We identify $\CC$ with a full subcategory of the $k$-linear category
$\CC^i$ via the Yoneda
functor.

We denote by $\CC[y]$ the category with the same objects as $\CC$ and with 
$\Hom_{\CC[y]}(c,c')=\Hom_{\CC}(c,c')\otimes_{\BZ}\BZ[y]$.

\smallskip
Given an additive category $\CA$, we denote by $\Comp(\CA)$ the category of
complexes of objects of $\CA$ and by $\Ho(\CA)$ its homotopy category.

\smallskip
Let $\CU$ be a monoidal category. We denote by $\CU^{\rev}$ the monoidal category
equal to $\CU$ as a category but with $u\otimes^{\rev}u'=u'\otimes u$.

Consider an action of $\CU$ on a category $\CA$ via the monoidal
functor $M:\CU\to\End(\CA)$. Let $\Sigma:\CU\to\CU$ be a monoidal functor. We denote
by $\Sigma^*\CA$ the category $\CA$ with the action of $\CU$ given by $M\Sigma$.

\smallskip
Given $E$ a functor, we write $E$ for $\id_E:E\to E$.

\section{$2$-Representations of $\Gsl_2^+$}
\subsection{Monoidal category}
\label{se:monoidalcat}
Let $\CU$ be the ($k$-enriched) monoidal category generated by an object $E$ and arrows
$x:E\to E$ and $\tau:E^2\to E^2$ modulo relations
\begin{equation}
	\label{eq:xtau}
\tau^2=0,\ (\tau E)\circ (E\tau)\circ (\tau E)=(E\tau)\circ (\tau E)\circ (E\tau),\
\tau\circ (xE)-(Ex)\circ\tau=1=(xE)\circ\tau-\tau\circ (Ex).
\end{equation}

There are isomorphisms of monoidal categories 
$$\omega:\CU\iso\CU^\opp,\ E\mapsto E,\ \tau\mapsto\tau,\ x\mapsto x$$
$$\text{and }\CU\iso\CU^\rev,\ E\mapsto E,\ \tau\mapsto\tau,\ x\mapsto -x.$$

\smallskip
Given $n\ge 0$, let $H_n$ be the $k$-algebra generated by $x_1,\ldots,x_n$,
$\tau_1,\ldots,\tau_{n-1}$ modulo relations
$$\tau_i^2=0, \tau_i\tau_j=\tau_j\tau_i \text{ if }|i-j|\neq 1 \text{ and }
\tau_i\tau_{i+1}\tau_i=\tau_{i+1}\tau_i\tau_{i+1}$$
$$x_rx_s=x_sx_r,\ \tau_ix_r=x_r\tau_i \text{ if }r\neq i,i+1 \text{ and }
\tau_i x_i-x_{i+1}\tau_i=1=x_i\tau_i-\tau_i x_{i+1}$$
for $1\le r,s\le n$ and $1\le i,j<n$.

We put $s_i=\tau_i(x_i-x_{i+1})-1$. There is an injective morphism of algebras
$k[x_1,\ldots,x_n]\rtimes\GS_n\to H_n,\ x_i\mapsto x_i, (i,i+1)\mapsto s_i$.

We identify $H_{n-1}$ with a subalgebra of $H_n$.
The left $H_{n-1}$-module $H_n$ is free with basis
	$(1,\tau_{n-1},\ldots,\tau_{n-1}\cdots\tau_1)$.

We put $H_{-1}=0$.

\smallskip
There is an isomorphism of algebras
$$H_n\iso \End_{\CU}(E^n)^\opp,\ x_i\mapsto E^{n-i}xE^{i-1},\ \tau_i\mapsto E^{n-i-1}\tau E^{i-1}.$$
It gives rise to a fully faithful functor
$$\bigoplus_n \Hom(E^n,-):\CU\to \bigoplus_{n\ge 0}H_n\mMod.$$

\subsection{Algebra embeddings}
%Given $m\ge 1$, we define an algebra embedding
%$$\iota_m:H_{n}[y]\to H_{m+n}:
%x_i\mapsto x_{m+i},\ \tau_i\mapsto\tau_{m+i},\ y\mapsto x_m.$$
%
%Note that $H_{m+n}\iota_m$ is a $(H_{m+n},H_n[y])$-bimodule. We extend it to a
%$(H_{m+n}[y],H_n[y])$ by letting the left action of $y$ to be the right action of $y$.
%
%We define an algebra embedding
%$$\iota_0:H_{n}[y]\to H_{n+1}:
%x_i\mapsto x_{i},\ \tau_i\mapsto\tau_{i},\ y\mapsto x_{n+1}.$$
%We extend the $(H_{n+1},H_n[y])$-bimodule structure of $H_{n+1}\iota_0$ to that of a
%$(H_{n+1}[y],H_n[y])$-bimodule by letting the left action of $y$ to be the right action of $y$.
%
%
%\medskip

Given $m,n\ge 0$, $r\in\{1,\ldots,m\}\bigcup \BZ_{>m+n}$ and $s\ge\sup(m+n,r)$,
we define an injective algebra morphism
$$\iota_{m,r}:H_n[y]\to H_s,\ x_i\mapsto x_{m+i},\ \tau_i\mapsto \tau_{m+i},\
y\mapsto x_r.$$
We put $\iota_0=\iota_{0,n+1}$ and $\iota_m=\iota_{m,m}$ for $m>0$.

\smallskip
Note that $H_s\iota_{m,r}$ is a $(H_s,H_n[y])$-bimodule. We extend it to a
$(H_s[y],H_n[y])$-bimodule by letting the left action of $y$ be the right action of $y$.

Similarly, we extend the structure of $(H_n[y],H_s)$-bimodule on 
$\iota_{m,r}H_s$ to a structure of $(H_n[y],H_s[y])$-bimodule by letting the right
action of $y$ be the left action of $y$.
%\smallskip
%Fix $n\ge 1$. We define the $(H_n[y],H_{n-1}[y])$-bimodule 
%$H'_n$. The underlying $(H_n,H_{n-1})$-bimodule is $H_n$. The left and right actions
%of $y$ are given by right multiplication by $x_n$.
%We put $H'_{-1}=H'_0=0$.
%
%We define the $(H_{n-1}[y],H_n[y])$-bimodule 
%$H''_n$. The underlying right $H_n$-module is $H_n$. The left and right actions of $y$ are
%given by left multiplication by $x_1$. 
%The left action of $a\in H_{n-1}$ is given by left multiplication by $\iota_1(a)$.
%We put $H''_{-1}=H''_0=0$.

\medskip
We define the morphism of algebras
$$\iota^y_m:H_n[y]\to H_{m+n}[y],\ x_i\mapsto x_{m+i},\ \tau_i\mapsto\tau_{m+i},\ y\mapsto y.$$

%$$\iota^\circ_m:({^0H}_{n})^\opp\to{^0H}_{m+n}:
%x_i\mapsto x_{m+i},\ \tau_i\mapsto\tau_{m+i}.$$

%$$\iota'_{m-1,m+n+1}:{^0H}_{m-1}[y]\mapsto{^0H}_{m+n+1}:
%x_i\mapsto x_{n+i+1},\ \tau_i\mapsto\tau_{n+i+1},\ y\mapsto %x_{m+n+1}.$$

\subsection{$\CU$-modules}
A $\CU$-module is a category $\CV$ endowed with an action of $\CU$, ie, a monoidal functor
$\CU\to\End(\CV)$. An action of $\CU$ on a category $\CV$
is the same data (equivalence of $2$-categories) as the data of an endofunctor $E$ of $\CV$
and of $x\in\End(E)$ and $\tau\in\End(E^2)$ satisfying (\ref{eq:xtau}).

Consider two $\CU$-modules $(\CV,E,x,\tau)$ and $(\CV',E',x',\tau')$. A morphism of
$\CU$-modules 
$$(\CV,E,x,\tau)\to (\CV',E',x',\tau')$$
is the data of a pair $(\Phi,\varphi)$ where
$\Phi:\CV\to\CV'$ is a functor and $\varphi:\Phi E\iso E'\Phi$ is an isomorphism such that
$\varphi\circ(\Phi x)=(x'\Phi)\circ\varphi$ and $(E'\varphi)\circ(\varphi E)\circ (\Phi\tau)=
(\tau'\Phi)\circ(E'\varphi)\circ(\varphi E)$.

\smallskip
A right $\CU$-module is defined to be a $\CU^{\mathrm{rev}}$-module.

\subsection{Regular $2$-representation}
The left and right actions of $\CU$ on itself by tensor product give rise to commuting left and
right actions of $\CU$ on the abelian category
$\CA=\bigoplus_{n\ge 0}H_n\mMod$. Let us describe those actions.

For the left (resp. right) action, $E$ acts as the direct sum of functors 
$H_{n+1}\iota_l\otimes_{H_n}-:H_n\mMod\to H_{n+1}\mMod$ with $l=0$ (resp. $l=1$).
The endomorphism
$x$ acts as right multiplication by $x_{n+1}$ (resp. $x_1$) on $H_{n+1}$. Via
the isomorphism $H_{n+2}\iota_l\otimes_{H_{n+1}}H_{n+1}\iso H_{n+2},\ a\otimes b\mapsto
a\iota_l(b)$, the endomorphism $\tau$ acts on $H_{n+2}$ as right multiplication by
$\tau_{n+1}$ (resp. $\tau_1$).

\subsection{Vector $2$-representation}
Let $\CL_r=k[y]\mMod$ for $r=0,1$. There is an action of $\CU$ on $\CL=\CL_0\oplus\CL_1$. The
functor $E$ acts as $\Id:\CL_0\to\CL_1$ and $x$ acts as multiplication by $y$.

\section{Description of the tensor product}
\subsection{Tensor product}
We consider the tensor product $\CL\dotimes\CA$ for the action of $\CU$ \cite{Rou3}. This
is an abelian category with an $\infty$-categorical action of $\CU$ on $D^b(\CL\dotimes\CA)$.
We will only use here the classical action.

We will make a key use of the fact that the canonical morphism of $\CU$-modules
$\CL\to\CL_0$ (where $E$ acts by $0$ on $\CL_0$) induces a
morphism of $\CU$-modules
\begin{equation}
	\label{eq:extension}
D^b(\CL\dotimes\CA)\to
D^b(\CL_0\dotimes\CA).
\end{equation}

We will describe the category $\CB=\CL\dotimes\CA$ and the action of $E$ following
\cite{Rou3}. The actions of $x$ and $\tau$ will be described later using (\ref{eq:extension}).

	\begin{rem}
	The general tensor product construction for $\CU$ provides a dg-category. Because
	the $2$-representation of $\Gsl_2^{>0}$ extends to $\Gsl_2^{\ge 0}$, that dg-category is
	canonically
	equivalent to categories of complexes for an algebra with zero differential, cf \cite{Rou3}.
	This is the model we will use for $\CB$.
	\end{rem}
\subsection{Underlying category}

Given $n\ge 0$, let
$\CB_n$ be the category with
objects $\left[{\xy (0,5)*{M}, (0,-5)*{N},\ar^{\gamma}(0,-3)*{};(0,3)*{}, \endxy}\right]$, where $M$ is
an $H_n[y]$-module, $N$ an $H_{n-1}[y]$-module, and $\gamma:N\to M$ is
	a morphism of $H_{n-1}[y]$-modules  such that  $(y-x_n)\gamma(m)=0$ for all $m\in N$.
	We define
	$$\Hom_{\CB_n}(\left[{\xy (0,5)*{M}, (0,-5)*{N},\ar^{\gamma}(0,-3)*{};(0,3)*{}, \endxy}\right],\left[{\xy (0,5)*{M'}, (0,-5)*{N'},\ar^{\gamma'}(0,-3)*{};(0,3)*{}, \endxy}\right])$$
		to be the set of pairs $(f,g)$ in
		$\Hom_{H_n[y]}(M,M')\oplus\Hom_{H_{n-1}[y]}(N,N')$ such that
		$f\circ\gamma=\gamma'\circ g$.

		Let
$$P_n^+=\left[{\xy (0,5)*{H_n[y]},(0,-5)*{0}, \ar(0,-3)*{};(0,3)*{},\endxy}\right]\text{ and }
P_n^-=\left[{\xy (0,5)*{H_n\iota_0}, (0,-5)*{H_{n-1}[y]}\ar^{\can}(0,-3)*{};(0,3)*{},\endxy}\right].$$
We have
$$\Hom(P_n^+,\left[{\xy (0,5)*{M}, (0,-5)*{N},\ar^{\gamma}(0,-3)*{};(0,3)*{}, \endxy}\right])
\iso M,\ (f,g)\mapsto f(1)$$
$$\Hom(P_n^-,\left[{\xy (0,5)*{M}, (0,-5)*{N},\ar^{\gamma}(0,-3)*{};(0,3)*{}, \endxy}\right])
\iso N,\ (f,g)\mapsto g(1).$$
So $P_n^+\oplus P_n^-$ is a progenerator of $\CB_n$.
Right multiplication provides isomorphisms of algebras
$$H_n[y]\iso\End(P_n^+)^\opp \text{ and }H_{n-1}[y]\iso\End(P_n^-)^{\opp}.$$

		\smallskip
We put $\CB=\bigoplus_{n\ge 0}\CB_n$.

\subsection{Left action}
We describe here the functor providing the action of the generator $E$ of $\CU$ 
 on $\CB$, following \cite{Rou3}.

\begin{prop}
	\label{pr:leftE}
We have a functor $E:\CB_n\to \Comp^b(\CB_{n+1})$ given by
$$E(\left[{\xy (0,5)*{M}, (0,-5)*{N},\ar^{\gamma}(0,-3)*{};(0,3)*{}, \endxy}\right])=
\left[{\xy (0,5)*{0\to H_{n+1}[y]\otimes_{H_n[y]} M}, (55,5)*{H_{n+1}\iota_0\otimes_{H_n[y]}M\to 0},
(0,-5)*{0\to H_n[y]\otimes_{H_{n-1}[y]}N},(55,-5)*{M\to 0}, 
\ar_{a\otimes m\mapsto a\gamma(m)}(30,-5)*{};(41,-5)*{}, \ar^{\can}(23,5)*{};(35,5)*{},
\ar^{a\otimes m\mapsto a(x_n-y)\tau_n\otimes \gamma(m)}(0,-3)*{};(0,3)*{},
\ar_{m\mapsto 1\otimes m}(55,-3)*{};(55,3)*{}, \endxy}\right]$$
where the non-zero terms of the complexes are in degrees $0$ and $1$
and
$$E(\left[\begin{matrix}f\\g\end{matrix}\right])=
\left[\begin{matrix}
	1\otimes f & 1\otimes f \\
	1\otimes g & f
\end{matrix}\right].
	$$
	The functor $E$ is exact.
\end{prop}

\begin{proof}
	Note first that the map
	$$u:H_n[y]\otimes_{H_{n-1}[y]}N\to H_{n+1}[y]\otimes_{H_n[y]} M,\
	a\otimes m\mapsto a(x_n-y)\tau_n\otimes \gamma(m)$$
	is well defined and it is a morphisn of $H_n[y]$-modules.
	We have
	$$(y-x_{n+1})u(a\otimes m)=a(x_n-y)(y-x_{n+1})\tau_n\otimes\gamma(m)=
a\tau_n(y-x_{n+1})\otimes (x_n-y)\gamma(m)=0.$$

Let us show that $E$ defines a complex.
	Given $a\in H_n[y]$ and $m\in N$, we have
	$$u(a\otimes m)=
	a(\tau_n(x_{n+1}-y)+1)\otimes\gamma(m),$$
	and this has image $a\otimes\gamma(m)$ in $H_{n+1}\iota_0\otimes_{H_n[y]}M$. It
	follows that $E$ defines a complex of objects of $\CB_{n+1}$.

Since $H_{n+1}$ is projective as a right $H_n[x_{n+1}]$-module, it follows that
$H_{n+1}\iota_0$ is projective as a right $H_n[y]$-module. We deduce that $E$ is exact.
\end{proof}

We extend $E$ to an exact
functor $E:\Comp^b(\CB_n)\to \Comp^b(\CB_{n+1})$ by taking total complexes
and we still denote by $E$ the endofunctor of $\Comp^b(\CB)$ obtained as the sum of the functors
$E$ for each $n$.

Note that the action of $E$ on $\Comp^b(\CB)$ does not extend to an action
of $\CU$. It only gives rise to an action of $\CU$ on $D^b(\CB)$.

\subsection{Compatibility with the filtration of $\CL$}
\label{se:filt}
We have a canonical morphism of $\CU$-modules $\CL\to \CL_0$ given by
the projection. It induces
a morphism of $\CU$-modules $D^b({\protect\CL\dotimes\CA})\to D^b(\CL_0\dotimes\CA)$.
Via the canonical equivalence
$\CA[y]\iso\CL_0\dotimes\CA$, we obtain
a $\CU$-module functor $ D^b(\CL\dotimes\CA)\to D^b(\CA[y])$. Let us
describe this explicitely.

\medskip
Consider $\Upsilon,\Upsilon_-:\CB\to\CA[y]$ given by
$$\Upsilon(\left[{\xy (0,5)*{M}, (0,-5)*{N},\ar^{\gamma}(0,-3)*{};(0,3)*{}, \endxy}\right])=M
\text{ and }
\Upsilon_-(\left[{\xy (0,5)*{M}, (0,-5)*{N},\ar^{\gamma}(0,-3)*{};(0,3)*{}, \endxy}\right])=N.
$$

There is a quasi-isomorphism $\varphi_+$
$${\xy
(0,10)*{H_{n+1}[y]\otimes_{H_n[y]}\Upsilon(\left[{\xy (0,5)*{M}, (0,-5)*{N},\ar^{\gamma}(0,-3)*{};(0,3)*{}, \endxy}\right])=
H_{n+1}[y]\otimes_{H_n[y]}M}, (60,10)*{0}, 
(10,-10)*{\Upsilon E(\left[{\xy (0,5)*{M}, (0,-5)*{N},\ar^{\gamma}(0,-3)*{};(0,3)*{}, \endxy}\right])
=H_{n+1}[y]\otimes_{H_n[y]} M}, (70,-10)*{H_{n+1}\iota_0\otimes_{H_n[y]}M},
\ar(42,10)*{};(52,10)*{},
\ar_{\can}(46,-10)*{};(52,-10)*{},
\ar^{a\otimes m\mapsto a(x_{n+1}-y)\otimes m}(26,6)*{};(26,-6)*{},
\ar_{\varphi_+}(-12,5)*{};(-12,-4)*{},
\endxy}$$
because 
$$0\to H_{n+1}[y]\xrightarrow{a\mapsto a(x_{n+1}-y)}H_{n+1}[y]
\xrightarrow{\can} H_{n+1}\iota_0\to 0$$
is a split exact sequence of projective right $H_n[y]$-modules.

\smallskip
This gives rise to a morphism of $\CU$-modules
	$$D^b(\CB) \xrightarrow{(\Upsilon,\varphi_+^{-1})} D^b(\CA[y]).$$

\subsection{Right action}
\label{se:rightaction}

		The right action of $\CU$ on itself by tensor product induces a right
		action of $\CU$ on $\CL\dotimes\CU$, hence on $\CB$.
It is given by
$$(\left[{\xy (0,5)*{M}, (0,-5)*{N},\ar^{\gamma}(0,-3)*{};(0,3)*{}, \endxy}\right])E=
\left[{\xy (0,5)*{H_{n+1}[y]\iota_1^y\otimes_{H_n[y]}M}, (0,-5)*{H_{n}[y]\iota_1^y
\otimes_{H_{n-1}[y]}N},\ar^{1\otimes \gamma}(0,-3)*{};(0,3)*{}, \endxy}\right].$$
	The endomorphism $x$ acts by 
	$$\left[\begin{matrix} a\otimes m \\ a'\otimes m'\end{matrix}\right]\mapsto
	\left[\begin{matrix} ax_1\otimes m \\ a'x_1\otimes m'\end{matrix}\right].$$
		Multiplication gives an isomorphism
$$(\left[{\xy (0,5)*{M}, (0,-5)*{N},\ar^{\gamma}(0,-3)*{};(0,3)*{}, \endxy}\right])
E^2\iso
\left[{\xy (0,5)*{H_{n+2}[y]\iota_2^y\otimes_{H_n[y]}M}, (0,-5)*{H_{n+1}[y]\iota_2^y
\otimes_{H_{n-1}[y]}N},\ar^{1\otimes \gamma}(0,-3)*{};(0,3)*{}, \endxy}\right]$$
	and $\tau$ acts by
$$\left[\begin{matrix} a\otimes m \\ a'\otimes m'\end{matrix}\right]\mapsto
	\left[\begin{matrix} a\tau_1\otimes m \\ a'\tau_1\otimes m'\end{matrix}\right].$$

			\smallskip
		Note that all the functors constructed in \S\ref{se:filt}
		are compatible with the right action
		of $\CU$.

\subsection{Commutation of the actions}
There is an isomorphism from
$$\bigl(E(\left[{\xy (0,5)*{M}, (0,-5)*{N},\ar^{\gamma}(0,-3)*{};(0,3)*{}, \endxy}\right])\bigr)E=$$
$$=
\left[{\xy (0,8)*{H_{n+2}[y]\iota_1^y\otimes_{H_{n+1}[y]}H_{n+1}[y]\otimes_{H_n[y]} M},
(85,8)*{H_{n+2}[y]\iota_1^y\otimes_{H_{n+1}[y]}H_{n+1}\iota_0\otimes_{H_n[y]}M},
(0,-8)*{H_{n+1}[y]\iota_1^y\otimes_{H_{n}[y]}H_n[y]\otimes_{H_{n-1}[y]}N},
(85,-8)*{H_{n+1}[y]\iota_1^y\otimes_{H_{n}[y]}M}, 
\ar_{a\otimes b\otimes m\mapsto a\otimes b\gamma(m)}(40,-8)*{};(61,-8)*{}, \ar^{\can}(35,8)*{};(51,8)*{},
\ar^{a\otimes b\otimes m\mapsto a\otimes b(x_n-y)\tau_n\otimes \gamma(m)}(0,-6)*{};(0,6)*{},
\ar_{a\otimes m\mapsto a\otimes 1\otimes m}(75,-6)*{};(75,6)*{}, \endxy}\right]$$
to
$$L=\left[{\xy (0,8)*{H_{n+2}[y]\iota_1^y\otimes_{H_n[y]} M},
(75,8)*{H_{n+2}\iota_{1,n+2}\otimes_{H_n[y]}M},
(0,-8)*{H_{n+1}[y]\iota_1^y\otimes_{H_{n-1}[y]}N},
(75,-8)*{H_{n+1}[y]\iota_1^y\otimes_{H_{n}[y]}M}, 
\ar_{a\otimes m\mapsto a\otimes \gamma(m)}(30,-8)*{};(51,-8)*{}, \ar^{\can}(35,8)*{};(47,8)*{},
\ar^{a\otimes m\mapsto a(x_{n+1}-y)\tau_{n+1}\otimes \gamma(m)}(0,-5)*{};(0,5)*{},
\ar_{a\otimes m\mapsto a\otimes m}(75,-5)*{};(75,5)*{}, \endxy}\right]$$
given by
$$\left[\begin{matrix}a_1\otimes b_1\otimes m_1 & a_2\otimes b_2\otimes m_2 \\
a_3\otimes b_3\otimes m_3 & a_4\otimes m_4 \end{matrix}\right]\mapsto
\left[\begin{matrix}a_1\iota_1^y(b_1)\otimes m_1 & a_2\iota_1^y(b_2)\otimes m_2 \\
a_3\iota_1^y(b_3)\otimes m_3 & a_4\otimes m_4 \end{matrix}\right]$$
%for $m_1,m_2,m_4\in M_n$, $m_3\in M_{n-1}$, $a_1,a_2\in H_{n+2}$, $a_3,a_4,b_1,b_2\in H_{n+1}$ and
%$b_3\in H_n$.

\smallskip
We have also an isomorphism from
$$E\bigl((\left[{\xy (0,5)*{M}, (0,-5)*{N},\ar^{\gamma}(0,-3)*{};(0,3)*{}, \endxy}\right])E\bigr)=$$
$$=
\left[{\xy (0,8)*{H_{n+2}[y]\otimes_{H_{n+1}[y]}H_{n+1}[y]\iota_1^y\otimes_{H_n[y]} M},
(85,8)*{H_{n+2}\iota_0\otimes_{H_{n+1}[y]}H_{n+1}[y]\iota_1^y\otimes_{H_n[y]}M},
(0,-8)*{H_{n+1}[y]\otimes_{H_{n}[y]}H_n[y]\iota_1^y\otimes_{H_{n-1}[y]}N},
(85,-8)*{H_{n+1}[y]\iota_1^y\otimes_{H_{n}[y]}M}, 
\ar_{a\otimes b\otimes m\mapsto ab\otimes \gamma(m)}(35,-8)*{};(61,-8)*{}, \ar^{\can}(36,8)*{};(48,8)*{},
\ar^{a\otimes b\otimes m\mapsto a(x_{n+1}-y)\tau_{n+1}\otimes b\otimes \gamma(m)}(0,-5)*{};(0,5)*{},
\ar_{a\otimes m\mapsto 1\otimes a\otimes m}(75,-5)*{};(75,5)*{}, \endxy}\right]$$
to $L$
given by
$$\left[\begin{matrix}a_1\otimes b_1\otimes m_1 & a_2\otimes b_2\otimes m_2 \\
a_3\otimes b_3\otimes m_3 & a_4\otimes m_4 \end{matrix}\right]\mapsto
\left[\begin{matrix}a_1b_1\otimes m_1 & a_2b_2\otimes m_2 \\
a_3b_3\otimes m_3 & a_4\otimes m_4 \end{matrix}\right].$$
Composing the second isomorphism with the inverse of the first one, we obtain an isomorphism 
of functors
$$E\bigl((?)E\bigr)\iso \bigl(E(?)\bigr)E.$$
%This isomorphism commutes with the actions of $X$ and $\tau$.

\section{New $t$-structure}

\subsection{Hecke calculations}
Fix $n\ge 0$. We put $\iota_1H_0=0$.
We define 
$$\tilde{\Delta}_n=\sum_{1\le r\le n}\tau_r\cdots\tau_{n-1}\otimes
	\tau_1\cdots\tau_{r-1}\in H_n[y]\otimes_{H_{n-1}[y]}\iota_1H_n$$
and we denote by $\Delta_n$ the image of $\tilde{\Delta}_n$ in 
$H_n\iota_0\otimes_{H_{n-1}[y]}\iota_1H_n$.

	\begin{lemma}
		\label{le:Deltacommutes}
		We have $a\Delta_n=\Delta_n a$ for all $a\in H_n$.
	\end{lemma}

	\begin{proof}
		Given $i\in\{1,\ldots,n-1\}$, we have
$$\tau_i\Delta_n=\sum_{1\le r<i}\tau_r\cdots\tau_{n-1}\tau_{i-1}\otimes
		\tau_1\cdots\tau_{r-1}+0+\tau_i\cdots\tau_{n-1}\otimes\tau_1\cdots\tau_i+
		\sum_{i+1<r\le n}\tau_r\cdots\tau_{n-1}\tau_i\otimes \tau_1\cdots\tau_{r-1}
		$$
		$$=\sum_{1\le r<i}\tau_r\cdots\tau_{n-1}\otimes\tau_i
		\tau_1\cdots\tau_{r-1}+
\tau_i\cdots\tau_{n-1}\otimes\tau_1\cdots\tau_i+
		\sum_{i+1<r\le n}\tau_r\cdots\tau_{n-1}\otimes\tau_{i+1} \tau_1\cdots\tau_{r-1}
		$$
$$=\sum_{1\le r<i}\tau_r\cdots\tau_{n-1}\otimes
		\tau_1\cdots\tau_{r-1}\tau_i+
\tau_i\cdots\tau_{n-1}\otimes\tau_1\cdots\tau_i+
		\sum_{i+1<r\le n}\tau_r\cdots\tau_{n-1}\otimes\tau_1\cdots\tau_{r-1}\tau_i
		$$
		$$=\Delta_n \tau_i.$$

		\medskip
		Consider now $i\in\{1,\ldots,n\}$.
		Given $n\ge d\ge j\ge 1$, we have
		$$x_j\tau_j\cdots\tau_{d-1}=\tau_j\cdots\tau_{d-1}x_d+
		\sum_{j\le s\le d-1}\tau_j\cdots\tau_{s-1}\tau_{s+1}\cdots\tau_{d-1},$$
		hence
		$$x_i\Delta_n=\sum_{1\le r<i}\tau_r\cdots\tau_{n-1}x_{i-1}\otimes
		\tau_1\cdots\tau_{r-1}-
		\sum_{1\le r<i}\tau_r\cdots\tau_{i-2}\tau_i\cdots\tau_{n-1}\otimes
		\tau_1\cdots\tau_{r-1}$$
		$$+\tau_i\cdots\tau_{n-1}x_n\otimes\tau_1\cdots\tau_{i-1}+
		\sum_{s=i}^{n-1}\tau_i\cdots\tau_{s-1}\tau_{s+1}\cdots\tau_{n-1}\otimes
		\tau_1\cdots\tau_{i-1}+\sum_{r=i+1}^n \tau_r\cdots\tau_{n-1}x_i
		\otimes \tau_1\cdots \tau_{r-1}.$$

	We have
		$$\sum_{r=1}^{i-1}\tau_r\cdots\tau_{n-1}x_{i-1}\otimes
		\tau_1\cdots\tau_{r-1}=\sum_{r=1}^{i-1}\tau_r\cdots\tau_{n-1}\otimes x_i
		\tau_1\cdots\tau_{r-1}=\sum_{r=1}^{i-1}\tau_r\cdots\tau_{n-1}\otimes
		\tau_1\cdots\tau_{r-1}x_i,$$
$$\sum_{r=i+1}^n \tau_r\cdots\tau_{n-1}x_i \otimes \tau_1\cdots \tau_{r-1}=
\sum_{r=i+1}^n \tau_r\cdots\tau_{n-1} \otimes x_{i+1}\tau_1\cdots \tau_{r-1}$$
		$$=\sum_{r=i+1}^n \tau_r\cdots\tau_{n-1} \otimes \tau_1\cdots \tau_{r-1}x_i-
		\sum_{r=i+1}^n \tau_r\cdots\tau_{n-1} \otimes \tau_1\cdots\tau_{i-1}
		\tau_{i+1}\cdots\tau_{r-1}
		$$
		and
		$$\tau_i\cdots\tau_{n-1}x_n\otimes\tau_1\cdots\tau_{i-1}=
		\tau_i\cdots\tau_{n-1}\otimes x_1\tau_1\cdots\tau_{i-1}$$
		$$=
		\tau_i\cdots\tau_{n-1}\otimes \tau_1\cdots\tau_{i-1}x_i+
		\sum_{s=1}^{i-1}\tau_i\cdots\tau_{n-1}\otimes \tau_1\cdots
		\tau_{s-1}\tau_{s+1}\cdots\tau_{i-1}.$$

We have
$$\sum_{s=i}^{n-1}\tau_i\cdots\tau_{s-1}\tau_{s+1}\cdots\tau_{n-1}\otimes
		\tau_1\cdots\tau_{i-1}=
		\sum_{s=i}^{n-1}\tau_{s+1}\cdots\tau_{n-1}\tau_i\cdots\tau_{s-1}
		\otimes
		\tau_1\cdots\tau_{i-1}$$
		$$=\sum_{s=i}^{n-1}\tau_{s+1}\cdots\tau_{n-1}\otimes\tau_{i+1}\cdots\tau_s
\tau_1\cdots\tau_{i-1}$$
		$$=\sum_{r=i+1}^{n}\tau_{r}\cdots\tau_{n-1}
		\tau_1\cdots\tau_{i-1}\otimes\tau_{i+1}\cdots\tau_{r-1}$$
and
$$\sum_{s=1}^{i-1}\tau_i\cdots\tau_{n-1}\otimes \tau_1\cdots
		\tau_{s-1}\tau_{s+1}\cdots\tau_{i-1}=
\sum_{s=1}^{i-1}\tau_i\cdots\tau_{n-1}\otimes \tau_{s+1}\cdots\tau_{i-1}
\tau_1\cdots \tau_{s-1} $$
$$=\sum_{s=1}^{i-1}\tau_i\cdots\tau_{n-1} \tau_s\cdots\tau_{i-2}\otimes\tau_1\cdots \tau_{s-1}$$
$$=\sum_{r=1}^{i-1}\tau_r\cdots\tau_{i-2}\tau_i\cdots\tau_{n-1}\otimes\tau_1\cdots \tau_{r-1}.$$

We deduce that $x_i\Delta_n=\Delta_n x_i$
	\end{proof}

	We define a morphism of left $H_n[y]$-modules
\begin{align*}
	\nu_n:\iota_1 H_{n+1}&\to H_n[y]\oplus H_n[y]\otimes_{H_{n-1}[y]}
	\iota_1 H_n \\
	\iota_1(a)\tau_1\cdots\tau_{r}&\mapsto\begin{cases}
		(a,a\tilde{\Delta}_n) & \text{ if }r=n \\
		(0,a(y-x_n)\otimes\tau_1\cdots\tau_r) & \text{ if }r<n.
	\end{cases}
\end{align*}
where $a\in H_n[y]$ and $r\in\{0,\ldots,n\}$ and a morphism of left $H_{n+1}[y]$-modules
\begin{align*}
	\nu'_n:H_{n+1}\iota_0\otimes_{H_{n}[y]}\iota_1 H_{n+1}&\to 
	H_{n+1}\iota_0\oplus H_{n+1}\iota_{0,n}\otimes_{H_{n-1}[y]}\iota_1 H_n \\
	a\otimes \tau_1\cdots\tau_r&\mapsto
	\begin{cases}
		(a,-a\sum_{s=1}^n \tau_s\cdots\tau_{n-1}s_n\otimes\tau_1\cdots\tau_{s-1})
		&  \text{ if }r=n \\
		(0,a(x_n-x_{n+1})s_n\otimes\tau_1\cdots\tau_r) &
		\text{ if }r<n.
	\end{cases}
\end{align*}
where $a\in H_{n+1}$ and $r\in\{0,\ldots,n\}$.

\smallskip
Note that there is a commutative diagram
\begin{equation}
	\label{eq:commnu'}
	\xymatrix{
	H_{n+1}\iota_0\otimes_{H_n[y]}\iota_1 H_{n+1}\ar[r]^{\nu'_n} 
	\ar[d]_{1\otimes\nu_n}&
H_{n+1}\iota_0\oplus H_{n+1}\iota_{0,n}\otimes_{H_{n-1}[y]}\iota_1 H_n  \\
H_{n+1}\iota_0\otimes_{H_n[y]}H_n[y]\oplus H_{n+1}\iota_0\otimes_{H_n[y]}
H_n[y]\otimes_{H_{n-1}[y]}\iota_1 H_n\ar[r]^-\sim_-{\mathrm{mult}}
& H_{n+1}\iota_0\oplus H_{n+1}\iota_0\otimes_{H_{n-1}[y]}\iota_1 H_n
\ar[u]_{\tiny\left(\begin{matrix} \id & 0\\ 0 & a\otimes b\mapsto -as_n\otimes b\end{matrix}
	\right)}
}
\end{equation}

\begin{lemma}
	\label{le:exactseqforqiso}
There is an exact sequence of $H_n[y]$-modules
	$$0\to \iota_1 H_{n+1}\xrightarrow{\nu_n} H_n[y]\oplus H_n[y]\otimes_{H_{n-1}[y]}
	\iota_1 H_n\xrightarrow{(a\mapsto -a\Delta_n,\can)} H_n\iota_0
	\otimes_{H_{n-1}[y]}\iota_1 H_n\to 0$$
	and an exact sequence of $H_{n+1}[y]$-modules
$$0\to H_{n+1}\iota_0\otimes_{H_{n}[y]}\iota_1 H_{n+1}\xrightarrow{\nu'_n}
H_{n+1}\iota_0\oplus H_{n+1}\iota_{0,n}\otimes_{H_{n-1}[y]}\iota_1 H_n
	\xrightarrow{(a\mapsto -a\Delta_n,\can)} \bar{H}_{n+1}\otimes_{H_{n-1}[y]}\iota_1 H_n\to 0$$
where $\bar{H}_{n+1}=H_{n+1}\iota_0/H_{n+1}(x_{n+1}-x_n)$.
	\end{lemma}

	\begin{proof}
		There is an exact sequence of $(H_n[y],H_{n-1}[y])$-bimodules
		$$0\to H_n[y]\xrightarrow{a\mapsto a(y-x_n)}H_n[y]\xrightarrow{\can}H_n\iota_0
		\to 0.$$
		Tensoring by the free $H_{n-1}[y]$-module $\iota_1H_n$, we obtain an
		exact sequence of $H_n[y]$-modules
	$$0\to H_n[y]\otimes_{H_{n-1}[y]}\iota_1H_n
	\xrightarrow{a\otimes b\mapsto a(y-x_n)\otimes b}
H_n[y]\otimes_{H_{n-1}[y]}\iota_1H_n\xrightarrow{\can}H_n\iota_0\otimes_{H_{n-1}[y]}\iota_1H_n
		\to 0.$$
		That exact sequence is a subcomplex of the complex of $H_n[y]$-modules
		\begin{multline}
			\label{eq:bigcomplex}
		0\to H_n[y]\oplus H_n[y]\otimes_{H_{n-1}[y]}\iota_1H_n
			\xrightarrow{\tiny\left(\begin{matrix} \id & 0 \\ 
				a\mapsto a\tilde{\Delta}_n
 &a\otimes b\mapsto a(y-x_n)\otimes b
		\end{matrix}\right)}
			\\
			H_n[y]\oplus
		H_n[y]\otimes_{H_{n-1}[y]}\iota_1H_n\xrightarrow{(a\mapsto -a\Delta_n,\can)}H_n\iota_0\otimes_{H_{n-1}[y]}\iota_1H_n
		\to 0
		\end{multline}
		and the quotient is acyclic. We deduce that (\ref{eq:bigcomplex}) is acyclic.

		\smallskip
		We have an isomorphism of $H_n[y]$-modules
		\begin{align*}
			\iota_1 H_{n+1}&\iso H_n[y]\oplus H_n[y]\otimes_{H_{n-1}[y]}\iota_1H_n\\
			\iota_1(a)\tau_1\cdots\tau_r&\mapsto
			\begin{cases}
				(a,0) & \text{ if }r=n \\
				(0,a\otimes \tau_1\cdots\tau_r) & \text{ if }r<n
			\end{cases}
		\end{align*}
		where $a\in H_n[y]$.
The composition of this map with the map $\left(\begin{matrix} \id & 0 \\ a\mapsto
	a\tilde{\Delta}_n
 &a\otimes b\mapsto a(y-x_n)\otimes b
		\end{matrix}\right)$ of (\ref{eq:bigcomplex}) is equal to $\nu_n$. We deduce
		that the first sequence of the lemma is exact.

		\medskip
		Tensoring the first exact sequence of the lemma with $H_{n+1}\iota_0$ gives
		an exact sequence of $H_{n+1}[y]$-modules
$$0\to H_{n+1}\iota_0\otimes_{H_{n}[y]}\iota_1 H_{n+1}\xrightarrow{\nu''_n}
H_{n+1}\iota_0\oplus H_{n+1}\iota_0\otimes_{H_{n-1}[y]}\iota_1 H_n
	\xrightarrow{(a\mapsto -a\Delta_n,\can)} \bar{H}_{n+1}\otimes_{H_{n-1}[y]}\iota_1 H_n\to 0$$
	where 
		$$\nu''_n(a\otimes \tau_1\cdots\tau_r)=
	\begin{cases}
		(a,\sum_{s=1}^n \tau_s\cdots\tau_{n-1}\otimes\tau_1\cdots\tau_{s-1})
		&  \text{ if }r=n \\
		(0,a(y-x_n)\otimes\tau_1\cdots\tau_r) &
		\text{ if }r<n.
	\end{cases}$$
for $a\in H_{n+1}$ and $r\in\{0,\ldots,n\}$.
		There is an isomorphism of $(H_{n+1}[y],H_{n-1}[y])$-bimodules
		$$f:H_{n+1}\iota_0\iso H_{n+1}\iota_{0,n},\
		a\mapsto -as_n.$$
		We have $\nu'_n=(\id,f\otimes 1)\circ\nu''_n$ (commutative diagram
		(\ref{eq:commnu'})) and we deduce that the second
		sequence of the lemma is exact.
	\end{proof}

	\begin{lemma}
		\label{le:Deltashiftn}
		The following diagram commutes
	$$\xymatrix{
		H_{n+1}[y]\ar[rrr]^-{\tiny\left(\begin{matrix} \can \\ a\mapsto a\Delta_n\end{matrix}\right)}
			&&& H_{n+1}\iota_0\oplus
	H_{n+1}\iota_{0,n}\otimes_{H_{n-1}[y]}\iota_1 H_n \\
	H_{n+1}[y]\ar[u]^{\id}\ar[rrr]_{a\mapsto a\Delta_{n+1}} &&& 
	H_{n+1}\iota_0\otimes_{H_n[y]}\iota_1 H_{n+1} \ar[u]_{\nu'_n}
}$$
\end{lemma}

\begin{proof}
	We have
	\begin{align*}
		\nu'_n(\Delta_{n+1})&=
	\sum_{r=1}^{n+1}\nu'_n(\tau_r\cdots\tau_n\otimes\tau_1\cdots\tau_{r-1})\\
		&=
	(1,-\sum_{s=1}^n \tau_s\cdots\tau_{n-1}s_n\otimes\tau_1\cdots\tau_{s-1})+
	\sum_{r=1}^n(0,\tau_r\cdots\tau_n(x_n-x_{n+1})s_n\otimes\tau_1\cdots\tau_{r-1})\\
		&=(1,
	\sum_{r=1}^n\tau_r\cdots\tau_{n-1}s_n^2\otimes\tau_1\cdots\tau_{r-1}).
	\end{align*}
\end{proof}

\begin{lemma}
	\label{le:formulaDeltato1}
	We have $s_{n-1}\cdots s_1(x_n-y)\cdots(x_2-y)\cdot\Delta_n=1\otimes 1$ in 
$H_n\iota_0\otimes_{H_{n-1}[y]}\iota_1 H_n$.
\end{lemma}

\begin{proof}
	We proceed by induction on $n$. The lemma is clear for $n=1$. Assume the lemma
	holds for $n$.

	By Lemma \ref{le:Deltashiftn}, we have
	$$\nu'_n(s_n\cdots s_1(x_{n+1}-y)\cdots(x_2-y)\cdot\Delta_{n+1})=
	(0,s_n\cdots s_1(x_{n+1}-y)\cdots(x_2-y)\cdot\Delta_n)=$$
$$= (0,s_n(x_{n+1}-y)s_{n-1}\cdots s_1(x_n-y)\cdots(x_2-y)\cdot\Delta_n)=
	(0,s_n(x_{n+1}-y)(1\otimes 1))=
	\nu'_n(1\otimes 1).$$
	The lemma follows now from the injectivity of $\nu'_n$ (Lemma \ref{le:exactseqforqiso}).
\end{proof}

\begin{lemma}
	\label{le:Deltagenerates}
	The element $\Delta_n$ generates $H_n\iota_0\otimes_{H_{n-1}[y]}\iota_1 H_n$ as a
	$H_n[y]$-module.
\end{lemma}

\begin{proof}
	Note that $1\otimes 1$ generates $H_n\iota_0\otimes_{H_{n-1}[y]}\iota_1 H_n$ as a
	$(H_n[y],H_n)$-bimodule. It follows from Lemma \ref{le:formulaDeltato1} that 
	$\Delta_n$ generates $H_n\iota_0\otimes_{H_{n-1}[y]}\iota_1 H_n$ as a
	$(H_n[y],H_n)$-bimodule, hence as a $H_n[y]$-module by Lemma \ref{le:Deltacommutes}.
\end{proof}

\subsection{Generators}

Given $n\ge 0$, we define a complex $Y_n$ of objects of $\CB_n$ in degrees $0$ and $1$
$$Y_n=\left[{\xy
(0,5)*{H_n[y]}, (38,5)*{H_n\iota_0\otimes_{H_{n-1}[y]}\iota_1 H_n}, (0,-5)*{0}, (38,-5)*{\iota_1H_n},
\ar^{a\mapsto a\Delta_n}(8,5)*{};(18,5)*{},
\ar(8,-5)*{};(18,-5)*{},
\ar_{a\mapsto 1\otimes a}(35,-2)*{};(35,2)*{},\endxy}\right].$$
Thanks to Lemma \ref{le:Deltacommutes}, right multiplication defines a morphism 
$$\gamma_n:H_n\to \End(Y_n)^\opp.$$

The canonical isomorphism $H_{n-1}[y]\otimes_{H_{n-1}[y]}\iota_1H_n\iso \iota_1 H_n$ induces an
isomorphism
$$(0\to P_n^+\xrightarrow{} P_n^-\otimes_{H_{n-1}[y]}\iota_1H_n\to 0)\iso Y_n.$$

\medskip
We have isomorphisms 
$$H_{n+1}[y]\otimes_{H_n[y]}H_n[y]\xrightarrow[\textrm{mult.}]{\sim} H_{n+1}[y],\
H_{n+1}\iota_0\otimes_{H_n[y]}H_n[y]\xrightarrow[\textrm{mult.}]{\sim}  H_{n+1}\iota_0,$$
$$H_{n+1}[y]\otimes_{H_n[y]}H_n\iota_0\otimes_{H_{n-1}[y]}\iota_1 H_n
\iso H_{n+1}\iota_{0,n}\otimes_{H_{n-1}[y]}\iota_1 H_n,\ ay^i\otimes b\otimes c\mapsto
abx_n^i\otimes c$$
for $a\in H_{n+1}$, $b,c\in H_n$ and $i\ge 0$ and
$$H_{n+1}\iota_0\otimes_{H_n[y]}H_n\iota_0\otimes_{H_{n-1}[y]}\iota_1 H_n\iso
\bar{H}_{n+1}\otimes_{H_{n-1}[y]}\iota_1H_n,\
a\otimes b\otimes c\mapsto ab\otimes c.$$

They give rise to an isomorphism

$$f'_n:E(Y_n)\iso$$
$$X_n=\left[{\xy
(0,7)*{H_{n+1}[y]},
(55,7)*{H_{n+1}\iota_0\oplus
H_{n+1}\iota_{0,n}\otimes_{H_{n-1}[y]}\iota_1 H_n},
(125,7)*{\bar{H}_{n+1}\otimes_{H_{n-1}[y]}\iota_1H_n},
(0,-7)*{0},
(55,-7)*{H_n[y]\oplus H_n[y]\otimes_{H_{n-1}[y]}\iota_1H_n},
(125,-7)*{H_n\iota_0\otimes_{H_{n-1}[y]}\iota_1H_n},
\ar^{\tiny\left(\begin{matrix} \can \\ a\mapsto a\Delta_n\end{matrix}\right)}(12,7)*{};(20,7)*{},
	\ar^{(a\mapsto -a\Delta_n,\can)}(87,7)*{};(100,7)*{},
\ar(8,-7)*{};(18,-7)*{},
\ar_{(a\mapsto -a\Delta_n, \can)}(85,-7)*{};(100,-7)*{},
\ar_{\tiny\left(\begin{matrix} \iota_0 & 0\\0 & a\otimes b\mapsto a(x_n-y)\tau_n\otimes
	b\end{matrix}\right)}(35,-5)*{};(35,5)*{},
\ar_{\can}(125,-5)*{};(125,5)*{},
\endxy}\right].$$

\begin{prop}
	\label{pr:YX}
	There is a morphism $f_n:Y_{n+1}\to X_n$ given by

$$\left[\begin{matrix}
	\id	& \nu'_n & 0 \\
0 &\nu_n & 0
\end{matrix}\right].$$
	The map $f_n$ is a quasi-isomorphism.
\end{prop}

\begin{proof}
	Since $(x_n-y)\tau_n=\tau_n(x_{n+1}-x_n)+1=-s_n$ in $H_{n+1}\iota_{0,n}$,
	the commutativity of the diagram  (\ref{eq:commnu'}) implies that the 
	following diagram commutes
	$$\xymatrix{
		H_{n+1}\iota_0\otimes_{H_n[y]}\iota_1 H_{n+1}\ar[r]^-{\nu'_n}
& H_{n+1}\iota_0\oplus
H_{n+1}\iota_{0,n}\otimes_{H_{n-1}[y]}\iota_1 H_n \\
	\iota_1 H_{n+1} \ar[r]_-{\nu_n} \ar[u]^{a\mapsto 1\otimes a}&
	H_n[y]\oplus H_n[y]\otimes_{H_{n-1}[y]}\iota_1H_n
	\ar[u]_-{\tiny\left(\begin{matrix} \iota_0 & 0\\0 & a\otimes b\mapsto a(x_n-y)\tau_n\otimes
	b\end{matrix}\right)}
}$$

	This shows that $(f_n)^1$ defines a morphism $(Y_n)^1\to (X_n)^1$.

	\smallskip

By Lemma \ref{le:exactseqforqiso}, the compositions
	$$\iota_1 H_n\xrightarrow{\nu_n} H_n[y]\oplus H_n[y]\otimes_{H_{n-1}[y]}\iota_1H_n
	\xrightarrow{(a\mapsto -a\Delta_n,\can)}
H_n\iota_0\otimes_{H_{n-1}[y]}\iota_1H_n$$
and 
	$$
	H_{n+1}\iota_0\otimes_{H_n[y]}\iota_1 H_{n+1}\xrightarrow{\nu'_n}
 H_{n+1}\iota_0\oplus
	H_{n+1}\iota_{0,n}\otimes_{H_{n-1}[y]}\iota_1 H_n \xrightarrow{(a\mapsto -a\Delta_n,\can)}
	\bar{H}_{n+1}\otimes_{H_{n-1}[y]}\iota_1H_n$$
	vanish.

	Together with Lemma \ref{le:Deltashiftn}, this shows that $f_n$ is a morphism of complexes.
	Lemma \ref{le:exactseqforqiso} shows that $f_n$ is a quasi-isomorphism.
\end{proof}

We put $g_n=f^{\prime -1}_n\circ f_n:Y_{n+1}\to E(Y_n)$. We define inductively
$h_n:Y_n\to E^{n}(Y_0)$ by $h_0=\id$ and $h_{n+1}=E(h_n)\circ g_n$.
As an immediate consequence of Proposition \ref{pr:YX}, we have the following corollary.

\begin{cor}
	\label{co:Ynqiso}
	The map $h_n:Y_n\to E^{n}(Y_0)$ is a quasi-isomorphism.
\end{cor}

\begin{lemma}
	\label{le:P+inj}
	Let $C$ and $C'$ be two complexes of objects of $\CB_n$ with
	\begin{itemize}
		\item $C^0$ and $C^{\prime 0}$ isomorphic to a direct sum of copies of $P^+_n$
		\item $C^1$ and $C^{\prime 1}$ isomorphic to a direct sum of copies of $P^-_n$
		\item $C^i=C^{\prime i}=0$ for $i\not\in\{0,1\}$.
	\end{itemize}
	The map $\Upsilon:\Hom_{D(\CB_n)}(C,C'))\to \Hom_{D(H_n[y])}(\Upsilon(C),\Upsilon(C'))$ is injective.
\end{lemma}

\begin{proof}
	The canonical map
	$\Hom_{\Ho(\CB)}(C,C')\to\Hom_{D(\CB)}(C,C')$ is an isomorphism.
	Note that since $\Hom(C^1,C^{\prime 0})=0$, the canonical map
	$\Hom_{\Comp(\CB)}(C,C')\to\Hom_{\Ho(\CB)}(C,C')$ is also an isomorphism.
 
	\smallskip
	Consider now $\alpha\in\Hom_{\Comp(\CB)}(C,C')$ such that $\Upsilon(\alpha)$ vanishes in
	$D(H_n[y])$. We have $H^0(\Upsilon(\alpha))=0$. Let
	$\kappa=\prod_{i=1}^n(y-x_i)\in k[x_1,\ldots,x_n]^{\GS_n}[y]\subset Z(H_n[y])$.
	We have $\kappa\cdot \Upsilon(P_n^-)=0$, hence
	$\kappa \Upsilon(C^0)\subset H^0(\Upsilon(C))$. We deduce that the restriction of
	$\Upsilon(\alpha^0)$ to $\kappa \Upsilon(C^0)$ vanishes.
	The algebra $H_n[y]$ is free as a 
	$k[x_1,\ldots,x_n,y]$-module, hence also as a $k[x_1,\ldots,x_n]^{\GS_n}[y]$-module. So
	the annihilator of $\kappa$ in $\Upsilon(C^{\prime 0})$ is $0$.
	It follows that $\Upsilon(\alpha)^0=0$.

The composition
	$\Upsilon(C^1)[-1]\xrightarrow{\Upsilon(\alpha^1)[-1]}\Upsilon(C^{\prime 1})[-1]
	\xrightarrow{\can}\Upsilon(C')$ is equal to the composition
	$\Upsilon(C^1)[-1]\xrightarrow{\can}\Upsilon(C)\xrightarrow{\Upsilon(\alpha)}\Upsilon(C')$, hence
	it vanishes in $D(H_n[y])$.
	There is an exact sequence
	$$\Hom_{H_n[y]}(\Upsilon(C^1),\Upsilon(C^{\prime 0}))\to
	\Hom_{H_n[y]}(\Upsilon(C^1),\Upsilon(C^{\prime 1}))\to
	\Hom_{D(H_n[y])}(\Upsilon(C^1),\Upsilon(C')[1]).$$
	Since $\kappa \Upsilon(C^1)=0$, it follows that $\Hom_{H_n[y]}(\Upsilon(C^1),\Upsilon(C^{\prime 0}))=0$.
	We deduce that $\Upsilon(\alpha^1)=0$, hence $\Upsilon(\alpha)=0$.

	Since the structure map $\Upsilon_-(C^{\prime 1})\to \Upsilon(C^{\prime 1})$
	is injective, it follows that $\Upsilon_-(\alpha)=0$. So $\alpha=0$.
\end{proof}

The structure of $\CU$-module on $D(\CB)$ provides a morphism of algebras 
$\gamma'_n:H_n\to\End_{D(\CB)}(E^{n}(Y_0))$.
The following proposition shows that the action of $H_n$ on $Y_n$ defined earlier is
compatible with $\gamma'_n$.

\begin{prop}
	\label{pr:actionYn}
We have a commutative diagram
	$$\xymatrix{
		& \End_{D(\CB)}(Y_n) \ar[dd]^{h_n}_\sim\\
		H_n\ar[ur]^-{\gamma_n}\ar[dr]_-{\gamma'_n}\\
	& \End_{D(\CB)}(E^{n}(Y_0))
	}$$
\end{prop}

\begin{proof}
	Multiplication provides an isomorphism $r$ from $(\Upsilon E^{n}(Y_0))^0$,
	the degree $0$ term of
	$\Upsilon E^{n}(Y_0)$, to $H_n[y]$. 
	By \S \ref{se:filt}, we have a quasi-isomorphism 
	$$\rho:H_n[y]\to \Upsilon E^{n}(Y_0),\ a\mapsto r^{-1}(a(x_n-y)\cdots (x_1-y))$$
and a commutative diagram
	$$\xymatrix{
		\End_{D(\CB)}(E^{n}(Y_0))\ar[r]^-{\Upsilon} &
		\End_{D(H_n[y])}(\Upsilon(E^{n}(Y_0))) \\
	H_n[y]\ar[u]^{\gamma'_n}\ar[r]_-{\text{right mult.}}& \End_{H_n[y]}(H_n[y])
	 \ar[u]^\sim_\rho
	}.$$

	There is a commutative diagram
	$$\xymatrix{
		& \Upsilon(Y_n)\ar[dd]^{\Upsilon(h_n)}\\
		H_n[y]\ar[ur]^-{a\mapsto a(x_n-y)\cdots(x_1-y)\ \ \ } \ar[dr]_-{\rho} \\
		& \Upsilon(E^{n}(Y_0))
	}$$
	We deduce that the following diagram commutes:
$$\xymatrix{
	& \End_{D(\CB)}(Y_n) \ar[r]^-{\Upsilon} & \End_{D(H_n[y])}(\Upsilon(Y_n))\ar[dd]^{\Upsilon(h_n)}_\sim \\
		H_n\ar[ur]^-{\gamma_n}\ar[dr]_-{\gamma'_n}\\
	& \End_{D(\CB)}(E^{n}(Y_0))\ar[r]_-{\Upsilon} & \End_{D(H_n[y])}(\Upsilon E^{n}(Y_0))
	}$$
The proposition follows from Lemma \ref{le:P+inj}.
\end{proof}

\subsection{$2$-representation on an additive subcategory}
Given $m,n\ge 0$, we put 
$$Y_{n,m}=\left[{\xy
(0,5)*{H_{m+n}[y]}, (43,5)*{H_{m+n}\iota_{m,m+n}\otimes_{H_{n-1}[y]}\iota_1 H_n}, (0,-5)*{0}, (43,-5)*{H_{m+n-1}[y]\iota_m^y\otimes_{H_{n-1}[y]}{\iota_1H_n}},
\ar^{d_{n,m}}(8,5)*{};(18,5)*{},
\ar(8,-5)*{};(18,-5)*{},
\ar_{\can}(43,-2)*{};(43,2)*{},\endxy}\right]$$
where $d_{n,m}(a)=a\sum_{r=1}^n\tau_{m+r}\cdots\tau_{m+n-1}\otimes \tau_1\cdots\tau_{r-1}$.
Note that $Y_n=Y_{n,0}$.

\smallskip
Right multiplication provides a morphism of algebras
$$H_n\to \End_{\Comp(\CB)}(Y_{n,m})^\opp.$$
$$H_m\to \End_{\Comp(\CB)}(Y_{n,m}),\
	h\mapsto \left(\left[\begin{matrix}a_0 & b_1\otimes c_1 \\
0 & b'_1\otimes c'_1\end{matrix}\right]\mapsto \left[\begin{matrix}a_0h & b_1h\otimes c_1 \\
0 & b'_1h\otimes c'_1\end{matrix}\right]\right).$$

Multiplication provides an isomorphism $(Y_n)E^m\iso Y_{n,m}$ that is equivariant for the 
action of $H_m$. Corollary \ref{co:Ynqiso} provides a quasi-isomorphism $Y_{n,m}\to E^n(Y_{0,m})$ 
compatible with the given maps $H_n\to \End_{D(\CB)}(Y_{n,m})^\opp$ and
$H_n\to \End_{D(\CB)}(E^n(Y_{0,m}))^\opp$ (Proposition \ref{pr:actionYn}).

\smallskip
We have an isomorphism given by multiplication
$$\bigl(0\to P_{m+n}^+\to
P_{m+n}^-\otimes_{H_{m+n-1}[y]}(H_{m+n-1}[y]\iota_m^y\otimes_{H_{n-1}[y]}\iota_1H_n)
\to 0\bigr)
\iso Y_{n,m}.$$

\begin{lemma}
	\label{le:Ydisjunct}
Given $m,n,m',n'\ge 0$, we have
	$\Hom_{D(\CB)}(Y_{n,m},Y_{n',m'}[i])=0$ for $i\neq 0$.
\end{lemma}

\begin{proof}
	Let $r=m+n$. We have a morphism 
	$$P_r^-\otimes_{H_{r-1}[y]}(H_{r-1}\iota_m^y\otimes_{H_{n-1}[y]}\iota_1H_n)\to
	P_r^-,\ a\otimes 1\otimes\tau_1\cdots\tau_{r-1}\mapsto
	\delta_{r,n}a \text{ for }1\le r\le n.$$
	Via the isomorphism above, it gives rise to a morphism
	$(Y_{n,m})^1\to P_r^-$ whose composition with $d_{m,n}$ is the canonical
	map $P_r^+\to P_r^-$.
	Since $\Hom(P_r^+,P_r^-)$ is generated by the canonical map
	as a module over $\End(P_r^-)$, 
	it follows that $\Hom_{\Ho(\CB)}(Y_{n,m},P_r^-)=0$, hence
	$\Hom_{\Ho(\CB)}(Y_{n,m},Y_{n',m'}[1])=0$,
	as $(Y_{n',m'})^1\simeq (P_{m'+n'}^-)^{\oplus n'}$.
	Since $\Hom(P_r^-,P_r^+)=0$, the lemma
	follows.
\end{proof}

\begin{lemma}
	\label{le:Ygenerate}
	The objects $Y_{0,m}$ and $Y_{1,m}$ for $m\ge 0$ generate $D^b(\CB)$ as a
	thick subcategory.
\end{lemma}

\begin{proof}
	The algebra $H_n[y]$ has finite global dimension since it is isomorphic
	to a matrix algebra over $k[x_1,\ldots,x_n]^{\GS_n}[y]$ (cf e.g.
	\cite[Proposition 2.21]{Rou2}).
	It follows that $\CB_n$ has finite global dimension since $\Hom(P_n^-,P_n^+)=0$
	and $\End(P_n^+)$ and $\End(P_n^-)$ have both finite global dimension and
	$P_n^+\oplus P_n^-$ is a progenerator for $\CB_n$.

	We have $Y_{0,m}=P_m^+$ and the cone of the canonical
	map $Y_{1,m-1}\to P_m^+$ is isomorphic to $P_m^-$. It follows that 
	$Y_{0,m}\oplus Y_{1,m-1}$ generates $D^b(\CB_n)$.
\end{proof}

Let $\CT$ be the full subcategory of $D(\CB)$ with objects the $Y_{n,m}$'s and $\CT'$ the 
full subcategory of $D(\CB)$ with objects those isomorphic to $Y_{n,m}$'s. 

Lemmas \ref{le:Ygenerate} and \ref{le:Ydisjunct} provide an equivalence 
$\Ho^b(\CT)\iso D^b(\CB)$.
Proposition \ref{pr:YX} shows that the action of $\CU$ on $D(\CB)$ restricts to an action of $\CT'$.
Also, the right action of $\CU$ on $D(\CB)$ restricts to $\CT'$. So, we obtain commuting left
and right actions of $\CU$ on $\CT$ and we have the following theorem.

\begin{thm}
	There is an equivalence of $(\CU,\CU)$-bimodules
$\Ho^b(\CT)\iso D^b(\CB)$.
\end{thm}

By Lemma \ref{le:P+inj}, the functor $\Upsilon:\CT\to D(\CA[y])$ is a faithful functor
of $(\CU,\CU)$-bimodules. Since $\Upsilon(Y_{0,0})=k[y]$, it follows that $\Upsilon(Y_{n,m})$ has
homology concentrated in degree $0$. Composing with the $H^0$ functor, we obtain a
faithful morphism of $(\CU,\CU)$-bimodules 
$$Q:\CT\to \CU[y],\ Y_{n,m}\mapsto E^{n+m}[y].$$

\section{Relation with Webster algebras}

\subsection{Webster category}
We follow \cite[\S 4.2]{We}.

Let $\tilde{\CW}$ be the free $(\CU,\CU)$-bimodule generated by
an object $\ast$ and by maps $\rho:\ast E\to E\ast$ and $\lambda:E\ast\to \ast E$.

Let $\CW$ be the  $(\CU,\CU)$-bimodule quotient of $\tilde{\CW}$ by the relations
$$\tau\ast\circ E\rho\circ \rho E=E\rho \circ \rho E\circ\ast\tau,\ 
\tau\ast\circ \lambda E\circ E\lambda=\lambda E \circ E\lambda\circ\ast\tau,$$
$$x\ast\circ \rho=\rho\circ \ast x,\ \ast x\circ \lambda=\lambda\circ x\ast,\ 
\lambda\circ\rho=\ast x,\ \rho\circ\lambda=x\ast$$
$$\rho E\circ \ast\tau\circ \lambda E-E\lambda \circ \tau \ast\circ E\rho=E\ast E.$$

\subsection{From Webster's category to the tensor product}
Let $\Sigma_y:\CU[y]\iso \CU[y]$ be the monoidal self-equivalence given by
$$\Sigma_y(E)=E,\ \Sigma_y(x)=x-y \text{ and }\Sigma_y(\tau)=\tau.$$

Let $\tilde{\Phi}:\tilde{\CW}[y]\to(\Sigma_y\otimes\Sigma_y)^*\CT$ be the (strict)
$(\CU,\CU)$-bimodules $k[y]$-enriched functor given by
$$\tilde{\Phi}(E^m\ast E^n)=Y_{m,n},\ \tilde{\Phi}(\rho)=\bar{\rho}
\text{ and }\tilde{\Phi}(\lambda)=\bar{\lambda}.$$

\begin{prop}
	The functor $\tilde{\Phi}$ factors through $\CW[y]$ and induces a morphism
	of $(\CU,\CU)$-bimodules $\Phi:\CW[y]\to(\Sigma_y\otimes\Sigma_y)^*\CT$.
\end{prop}
\begin{proof}
Note that
$$Y_{0,0}=Y_0=\left[\begin{matrix} k[y]\\0 \end{matrix}\right],\ 
Y_{1,0}=Y_1=\left[{\xy
(0,5)*{k[x_1,y]}, (25,5)*{k[x_1]}, (0,-5)*{0}, (25,-5)*{k[x_1]},
\ar^{\can}(8,5)*{};(18,5)*{},
\ar(8,-5)*{};(18,-5)*{},
\ar_{\id}(25,-2)*{};(25,2)*{},\endxy}\right]
\text{ and }
Y_{0,1}=\left[\begin{matrix} k[x_1,y]\\0 \end{matrix}\right].$$
	The action of $x$ on $Y_{1,0}=E(Y_0)$ (resp. on $Y_{0,1}=(Y_0)E$)
is given by multiplication by $x_1$.

We put
$$\bar{\lambda}=\left[\begin{matrix} \id \\ 0\end{matrix}\right]:Y_{1,0}\to Y_{0,1}
	\text{ and }
\bar{\rho}=\left[\begin{matrix} x_1-y \\ 0\end{matrix}\right]:Y_{0,1}\to Y_{1,0}.$$

	The quasi-isomorphism $E[y]=k[x_1,y]\to \Upsilon(Y_{1,0})$ is multiplication by $x_1-y$ in degree $0$.
	The isomorphism $E[y]=k[x_1,y]\to \Upsilon(Y_{0,1})$ is the identity.

It follows that $Q(\bar{\lambda})=x-y$ and $Q(\bar{\rho})=\id$.

\medskip
	We consider $\tilde{\Psi}=Q\circ\tilde{\Phi}:\tilde{\CW}[y]\to(\Sigma_y\otimes\Sigma_y)^*\CU[y]$.
	We have 
	$$\tilde{\Psi}(x*)\circ Q(\bar{\rho})=x-y=Q(\bar{\rho})\circ \tilde{\Psi}(*x)$$
	$$Q\bar{\lambda}\circ Q(\bar{\rho})=x-y=Q(*x)$$
	$$Q(\bar{\rho})\circ Q(\bar{\lambda})=x-y=Q(x*)$$
		$$\tilde{\Psi}(\tau*)\circ E Q(\bar{\rho})\circ Q(\bar{\rho})E=\tau
		= E Q(\bar{\rho})\circ Q(\bar{\rho})E\circ \tilde{\Psi}(*\tau) $$
		$$
	Q(\bar{\rho}) E\circ\tilde{\Psi}(\ast\tau)\circ Q(\bar{\lambda}) E-
	EQ(\bar{\lambda}) \circ \tilde{\Psi}(\tau \ast)\circ EQ(\bar{\rho})=
	\tau(xE-y)-(Ex-y)\tau=1.$$
	It follows that $\tilde{\Psi}$ factors through $\CW[y]$. Since $Q$ is faithful, the
	proposition follows.
\end{proof}

\begin{thm}
\label{th:Phiiso}
	The functor $\Phi$ is an isomorphism
	of $(\CU,\CU)$-bimodules
	$\CW[y]\iso(\Sigma_y\otimes\Sigma_y)^*\CT$.
\end{thm}

Theorem \ref{th:Phiiso} will be deduced from its graded version (Theorem \ref{th:Phiisograded}). Let
us first show the faithfulness of $\Phi$.

		\begin{lemma}
			\label{le:phifaithful}
			The functor $\Phi$ is faithful.
		\end{lemma}

		\begin{proof}
			We define a $k[u]$-linear functor
		$$\tilde{R}:\tilde{\CW}[y]\to(\Sigma_y\otimes\Sigma_y)^*\CU[y]$$
			of
			$(\CU,\CU)$-bimodules by
		$$\tR(\ast)=1,\ \tR(\lambda)=\tR(\rho)=E.$$
			The functor $\tR$ factors through $\CW[y]$ and induces a $k[y]$-linear functor
			$R:\CW[y]\to(\Sigma_y\otimes\Sigma_y)^*\CU[y]$ of 
		$(\CU,\CU)$-bimodules. It follows from \cite[Proposition 4.16]{We} that
			the functor $R$ is faithful.

			The composition $Q\circ\Phi:\CW[y]\to (\Sigma_y\otimes\Sigma_y)^*\CU[y]$ is
			a $k[y]$-linear functor
			of $(\CU,\CU)$-bimodules. The functors $Q\circ\Phi$ and $R$ take
			the same value on $\ast$ and on the generating arrows
			$\rho$ and $\lambda$. It follows that they are equal.

		$$\xymatrix{
			\CW[y]\ar[rr]^{R} \ar[dr]_\Phi && \CU[y] \\
			&\CT\ar[ur]_{Q}
	}$$

			Since $R$ is faithful, it follows that $\Phi$ is faithful.
		\end{proof}

	\subsection{Gradings}
		\subsubsection{Generalities}
Consider a category $\CC$ enriched in graded $k$-modules. We denote by $\CC\mgr$ the category
		enriched in $k$-modules with objects pairs $(c,n)$ where $c$ an object of $\CC$ and
		$n\in\BZ$ and with $\Hom_{\CC\mgr}((c,n),(c',n'))$ the space of homogeneous
		elements of degree $n'-n$ of $\Hom_{\CC}(c,c')$. This is a graded category, i.e.,
		a category endowed with an action of $\BZ$. We denote by $c\mapsto v^nc$ the action
		of $n\in\BZ$.

		For example, if $M$ is a graded $k$-vector space, then $(v^nM)_i=M_{i-n}$. When
		the homogenous components of $M$ are finite-dimensional, we put
		$\grdim(M)=\sum_{i\in\BZ}v^i\dim(M_i)$. We put $q=v^2$.

\smallskip
		Given $\CC$ a graded $k$-linear category, we equip $K_0(\CC)$ with a structure
		of $\BZ[v^{\pm 1}]$-module by $v[M]=[vM]$.
		When $\Hom$'s in $\CC$ are finite-dimensional, we define a bilinear form 
		$$K_0(\CC)\otimes_{\BZ}K_0(\CC)\to \BZ[v^{\pm 1}],\ 
		\langle [M],[N]\rangle=\grdim(\Hom(M,N)).$$

		\subsubsection{Gradings of $2$-representations}

	We enrich the monoidal category $\CU$ in graded $k$-vector spaces by setting $\deg(x)=2$ and
	$\deg(\tau)=-2$. Similarly, we define a grading on the algebra $H_n$ by setting
		$\deg(x_i)=2$ and $\deg(\tau_i)=-2$.

	We define $\deg(y)=2$.
We define $\CB'_n$ to be the category with
objects $\left[{\xy (0,5)*{M_n}, (0,-5)*{M_{n-1}},\ar^{\gamma}(0,-3)*{};(0,3)*{}, \endxy}\right]$,
		where $M_{r}$ is a graded $H_r[y]$-module and $\gamma$ is a morphism of graded 
	$H_{n-1}[y]$-modules $\gamma:M_{n-1}\to M_n$ such that  $(y-x_n)\gamma(m)=0$ for all $m\in M_{n-1}$.
		We define $\Hom_{\CB'_n}(M,N)$ to be the subspace of $\Hom_{\CB_n}(M,N)$ of
		graded morphisms. We define a graded version of the left action of $\CU$ by letting
		$E$ act by multiplying by $v$ the formula in Proposition \ref{pr:leftE}. The
		graded version of the left action of $\CU$ is obtained by using the formula in
		\S\ref{se:rightaction}. So, we have 
		a structure of $(\CU\mgr,\CU\mgr)$-bimodule on $\CB'$.

	We define a graded structure on $Y_{n,m}$ by
$$Y_{n,m}=\left[{\xy
	(0,5)*{H_{m+n}[y]}, (48,5)*{q^{n-1}H_{m+n}\iota_{m,m+n}\otimes_{H_{n-1}[y]}\iota_1 H_n}, (0,-5)*{0}, (48,-5)*{q^{n-1}H_{m+n-1}[y]\iota_m^y\otimes_{H_{n-1}[y]}{\iota_1H_n}},
\ar^{d_{m,n}}(8,5)*{};(18,5)*{},
\ar(8,-5)*{};(18,-5)*{},
\ar_{\can}(48,-2)*{};(48,2)*{},\endxy}\right]$$
	Corollary \ref{co:Ynqiso} provides a graded quasi-isomorphism $Y_{n,0}\to E^nY_{0,0}$.

	\smallskip
	There is a unique enrichment in graded $k$-vector spaces of $\CW$ that makes the
	structure of $(\CU,\CU)$-bimodule compatible with gradings and with 
	$\deg(\rho)=2$ and $\deg(\lambda)=0$. Note that while our choice of gradings differ
	from the one in \cite[Definition 4.4]{We}, the categories of graded modules are equivalent.

	\smallskip
The functor $\Phi$ gives rise to a morphism of graded $(\CU\mgr,\CU\mgr)$-bimodules
$\CW[y]\mgr\to(\Sigma_y\otimes\Sigma_y)^*\CT\mgr$.

\begin{lemma}
	\label{le:isometry}
	The functor $\Phi$ induces an isometry $K_0(\CW[y]\mgr)\iso K_0(\CT\mgr)$.
\end{lemma}

		\begin{proof}
		Let $\tilde{U}$ be the $\BQ(v)$-algebra generated by $e,f,k^{\pm 1}$ modulo relations
		$$ke=v^2ek,\ kf=v^{-2}fk\text{ and }
		ef-fe=\frac{k-k^{-1}}{v-v^{-1}}.$$

			We consider the coproduct on $\tilde{U}$ given by
		$$\Delta(e)=e\otimes k+1\otimes e,\ \Delta(f)=f\otimes 1+k^{-1}\otimes f
		\text{ and }\Delta(k)=k^{-1}.$$

\smallskip
		We put $[n]_v=\frac{v^n-v^{-1}}{v-v^{-1}}$ and $[n]_v!=[2]_v\cdots [n]_v$.

	Let $R=\BZ[v,v^{-1}]$ and let $U$ be the $R$-subalgebra of $\tilde{U}$ generated by
		the elements $e^{(n)}=\frac{e^n}{[n]_v!}$,
		$f^{(n)}=\frac{f^n}{[n]_v!}$ for $n\ge 1$ and by $k^{\pm 1}$.
		
		Let $U^+$ be the $R$-subalgebra of $U$  generated by
		the elements $e^{(n)}=\frac{e^n}{[n]_v!}$ for $n\ge 1$.
		There is an isomorphism of $R$-algebras
			$$U^+\iso K_0(\CU\ \mgr),\ e\mapsto [E].$$

		As a consequence, if $\CM$ is a graded category with a graded action of $\CU\mgr$, then
		$K_0(\CM)$ has a structure of $U^+$-module.

\smallskip
			Consider the $(n+1)$-dimensional representation $L(n)$ of $U$ with $R$-basis
		$(b_i)_{0\le i\le n}$ and
		$$f(b_i)=\delta_{i\neq 0} [n-i+1]_v b_{i-1},
		e(b_i)=\delta_{i\neq n} [i+1]_v b_{i+1} \text{ and }k(b_i)=v^{2i-n}.$$
%		The Shapovalov form on $L(n)$ is the $\BZ$-linear map $L(n)\times L(n)\to \BQ(v)$ given by
%		$$\langle v^ab_i,v^b b_j\rangle=
%		v^{b-a}\delta_{i,j}\frac{\prod_{r=1}^i(1-q^{n-r+1})}{\prod_{r=1}^i
%		(1-q^{r})}.$$

Let $V=L(1)$ and put $b_-=b_0$ and $b_+=b_1$.

\medskip
			We follow Webster, but we swap $e$ and $f$ and $k$ and $k^{-1}$.

		The $q$-Shapovalov form on $L(n)\otimes V$ is the $\BZ$-linear map 
$(L(n)\otimes V)\times (L(n)\otimes V)\to \BQ(v)$ defined by
		$$\langle v^a b_i\otimes e^c b_-,v^b b_j\otimes e^d b_-\rangle=
		v^{b-a}\delta_{i,j}\delta_{c,d}\frac{\prod_{r=1}^i(1-q^{n-r+1})}{\prod_{r=1}^i
		(1-q^{r})}.$$
for $a,b\in\BZ$, $i,j\in\{0,\ldots,n\}$ and $c,d\in\{0,1\}$.

\smallskip
We define a new form on $L(n)\otimes V$ by 
$$(u,u')=\langle u,u'+(v-v^{-1}) (e\otimes f)(u')\rangle.$$
So
$$(v^a b_i\otimes b_-,v^b b_j\otimes b_-)=
(v^a b_i\otimes b_+,v^b b_j\otimes b_+)=
v^{b-a}\delta_{i,j}\frac{\prod_{r=1}^i(1-q^{n-r+1})}{\prod_{r=1}^i (1-q^{r})},$$
$$(v^a b_i\otimes b_+,v^b b_j\otimes b_-)=0 \text{ and }
(v^a b_i\otimes b_-,v^b b_j\otimes b_+)=
-v^{b-a}\delta_{i-1,j}v^{-i}\frac{\prod_{r=1}^{i}(1-q^{n-r+1})}{\prod_{r=1}^{i-1} (1-q^{r})}$$
for $a,b\in\BZ$ and  $i,j\in\{0,\ldots,n\}$.

\smallskip
Taking the limit as  $n\to\infty$ (cf  \cite[Proof of Proposition 4.39]{We}),
we obtain a form on $U^+\otimes V$. This
 is the $\BZ$-linear map 
$(U^+\otimes V)\times (U^+\otimes V)\to \BQ(v)$ defined by
$$(v^a e^{(i)}\otimes b_-,v^b e^{(j)}\otimes b_-)=
(v^a e^{(i)}\otimes b_+,v^b e^{(j)}\otimes b_+)=
v^{b-a}\delta_{i,j}\frac{1}{\prod_{r=1}^i (1-q^{r})},$$
$$(v^a e^{(i)}\otimes b_+,v^b e^{(j)}\otimes b_-)=0 \text{ and }
(v^a e^{(i)}\otimes b_-,v^b e^{(j)}\otimes b_+)=
-v^{b-a}\delta_{i-1,j}v^{-i}\frac{1}{\prod_{r=1}^{i-1} (1-q^{r})}$$
for $a,b\in\BZ$ and  $i,j\in\BZ_{\ge 0}$.

\medskip
		We consider the $(U^+,U^+)$-bimodule $U^+\otimes V$ where
		the left action of $U^+$ is the diagonal action and the right action of
		$U^+$ is the action by right multiplication on $U^+$.
Webster \cite[Proposition 4.39]{We} shows there is an isomorphism of $(U^+,U^+)$-bimodules
$$\phi:U^+\otimes V\iso K_0((\CW\mgr)^i),\ 1\otimes b_-\mapsto [\ast]$$
that is compatible with the bilinear forms $(-,-)$. Note that Webster's isomorphism is
the composition of $\phi$ with the automorphism coming from the anti-involution of $\CW$ that
swaps $E^m\ast E^n$ with $E^n\ast E^m$, $\rho$ with $\lambda$, $\tau$ with $-\tau$ and that fixes $x$.

\medskip
		There is an isomorphism of $(U^+,U^+)$-bimodules
		$$\psi:U^+\otimes V\iso K_0(\CB\mgr),\ 1\otimes b_-\mapsto [P_0^+],$$
		where the $K_0$ is for $\CB\mgr$ as an abelian category.
		We have
		$$\psi(e^n\otimes b_-)=[P_n^+] \text{ and }\psi(e^n\otimes b_+)=-v^{-1}[P_{n+1}^-].$$

	We have 
	$$\grdim(H_n[y])=q^{-n(n-1)/2}\frac{(1-q)\cdots (1-q^n)}{(1-q)^{2n+1}},$$
hence	
	$$\langle P_i^+,P_j^+\rangle=\delta_{i,j}q^{-i(i-1)/2}\frac{(1-q)\cdots(1-q^i)}{(1-q)^{2i+1}}
	=\frac{1}{1-q}([i]_v!)^2 (e^{(i)}\otimes b_-,e^{(j)}\otimes b_-)$$
$$	\langle P_{i+1}^-,P_{j+1}^-\rangle=\delta_{i,j}q^{-i(i-1)/2}\frac{(1-q)\cdots(1-q^i)}{(1-q)^{2i+1}}
	=\frac{1}{1-q}([i]_v!)^2 (e^{(i)}\otimes b_+,e^{(j)}\otimes b_+)$$
$$	\langle P_i^+,P_j^-\rangle=\delta_{i,j}q^{-i(i-1)/2}\frac{(1-q)\cdots(1-q^i)}{(1-q)^{2i}}
	=\frac{-v}{1-q}[i]_v![i-1]_v! (e^{(i)}\otimes b_-,e^{(j-1)}\otimes b_+)$$
$$	\langle P_i^-,P_j^+\rangle=0$$

	We deduce that for all $w,w'\in U^+\otimes V$, we have 
	$(\psi(w),\psi(w'))=\frac{1}{1-q}(w,w')$.

	\smallskip
	The composition $\psi\circ\phi^{-1}$ is a morphism of $(U^+,U^+)$-bimodules sending
	$[\ast]$ to $[P_0^+]$. So, $\psi\circ\phi^{-1}$ and $[\Phi]$ agree on $[\ast]$. Since
	$1\otimes b_-$ generates $U^+\otimes V$ as a $(U^+,U^+)$-bimodule, we deduce that
	$\psi\circ\phi^{-1}=[\Phi]$ and the proposition follows.
\end{proof}

\begin{thm}
\label{th:Phiisograded}
	The functor $\Phi$ is an isomorphism
	of graded $(\CU\mgr,\CU\mgr)$-bimodules
	$\CW[y]\mgr\iso(\Sigma_y\otimes\Sigma_y)^*\CT\mgr$.
\end{thm}

\begin{proof}
	Lemma \ref{le:phifaithful} shows that $\Phi$ induces an injective map between graded $\Hom$-spaces.
	By Lemma \ref{le:isometry}, these $\Hom$ spaces have the same graded dimension. It follows that 
	$\Phi$ is fully faithful, hence it is an equivalence.
\end{proof}

\end{document}